\documentclass[11pt]{article}
\usepackage{amsmath,amsthm,amsfonts,amssymb,latexsym}
\usepackage{amsmath}
\usepackage{color}

\usepackage{a4wide}
\newtheorem{thm}{Theorem}[section]

\newtheorem{lem}[thm]{Lemma}

\newtheorem{theorem}{Theorem}[section]

\newtheorem{remark}[theorem]{Remark}

\def\k{\mathfrak k}

\def\l{\lambda}

\def\k{\mathfrak{k}}

\def\L{{\cal   L}}

\def\L1#1{L^1(#1)}

\def\lef({\left(}
\def\rig){\right)}

\numberwithin{equation}{section}

\begin{document}

\title{  Cohomology of  $\frak {sl}(2)$ acting
on the space of
$n$-ary differential operators on $\mathbb{R}$}

\label{firstpage}

\author{ Mabrouk Ben Ammar \and Rabeb Sidaoui \thanks{
 Universit\'e de Sfax, Facult\'e des Sciences, D\'epartement de Math\'ematiques, Laboratoire d'Alg\`ebre, G\'eom\'etrie et Th\'eorie
 Spectrale (AGTS) LR11ES53, BP 802, 3038 Sfax, Tunisie.
E.mail:  mabrouk.benammar@fss.rnu.tn}}

\maketitle

\begin{abstract}
 We consider the spaces $\mathcal{F}_\mu$ of polynomial  $\mu$-densities on the line as $\mathfrak{sl}(2)$-modules and then we compute the cohomological
 spaces $\mathrm{H}^1_\mathrm{diff}(\mathfrak{sl}(2), \mathcal{D}_{\bar{\lambda},\mu})$,   where $\mu\in \mathbb{R}$, $\bar{\lambda}=(\lambda_1,\dots,\lambda_n)
 \in\mathbb{R}^n$ and $\mathcal{D}_{\bar{\lambda},\mu}$ is the space of $n$-ary differential operators from $\mathcal{F}_{\lambda_1}\otimes\cdots\otimes
 \mathcal{F}_{\lambda_n}$ to $\mathcal{F}_\mu$.
\end{abstract}

\maketitle {\bf Mathematics Subject Classification} (2010). 17B56


\thispagestyle{empty}

\section{Introduction}

Consider the space of polynomial  $\mu$-densities:
\begin{equation*}
\mathcal{F}_\mu=\big\{ fdx^{\mu}, ~f\in \mathbb{R}[x]\big\},\quad \mu\in\mathbb{R}.
\end{equation*}
The Lie algebra $\mathrm{Vect}(\mathbb{R})$ of  polynomial vector fields $X_h=h{d\over dx}$, where $h\in \mathbb{R}[x]$, acts on $\mathcal{F}_\mu$ by the
{\it Lie derivative} $L^\mu$:
\begin{equation}\label{Lie1}X_h\cdot (fdx^{\mu}) =L_{X_h}^\mu(fdx^{\mu}):=(hf'+\mu h'f)dx^\mu.
\end{equation}

For $\bar{\lambda}=(\lambda_1,\dots,\lambda_n)\in\mathbb{R}^n$ and $\mu\in\mathbb{R}$ we denote by $\mathcal{D}_{\bar{\lambda},\mu}$ the space of $n$-ary
differential operators $A$ from $\mathcal{F}_{\l_1}\otimes\cdots\otimes\mathcal{F}_{\l_n}$ to $\mathcal{F}_\mu$. The Lie algebra $\mathrm{Vect}(\mathbb{R})$
acts on the space $\mathcal{D}_{\bar{\lambda},\mu}$ of these differential operators by:
\begin{equation}\label{Lieder2}
X_h\cdot A:= L_{X_h}^{\bar{\lambda},\mu}(A)=L_{X_h}^\mu\circ A-A\circ L_{X_h}^{\bar{\lambda}}
\end{equation}
where $L_{X_h}^{\bar{\lambda}}$ is the Lie derivative on $\mathcal{F}_{\lambda_1}\otimes\cdots\otimes\mathcal{F}_{\lambda_n}$  defined by the Leibnitz rule.
The spaces $\mathcal{F}_\mu$ and $\mathcal{D}_{\bar{\lambda},\mu}$ can be also viewed as $\mathfrak{sl}(2)$-modules, where $\mathfrak{sl}(2)$ is realized as
a subalgebra of $\mathrm{Vect}(\mathbb{R})$:
$$
\mathfrak{ sl}(2)=\mathrm{Span}(X_1,\,X_x,\,X_{x^2}).
$$

According to Nijenhuis-Richardson \cite{r1}, the space $\mathrm{H}^1\left(\mathfrak{g},\mathrm{End}(V)\right)$ classifies the infinitesimal deformations of a
 $\mathfrak{g}$-module $V$ and the obstructions to integrability of a given infinitesimal deformation of $V$ are elements of $\mathrm{H}^2\left(\mathfrak{g},
 \mathrm{End}(V)\right)$.

For $\bar{\lambda}\in\mathbb{R}$ the spaces $\mathrm{H}^1_\mathrm{diff}(\mathfrak{sl}(2), \mathcal{D}_{\bar{\lambda},\mu})$ are computed by Gargoubi \cite{g}
and Lecomte \cite{lec}, the spaces $\mathrm{H}^1_\mathrm{diff}(\mathrm{Vect}(\mathbb{R}),\mathfrak{sl}(2), \mathcal{D}_{\bar{\lambda},\mu})$ are computed by
Bouarroudj and Ovsienko \cite{bo} and the spaces $\mathrm{H}^1_\mathrm{diff}(\mathrm{Vect}(\mathbb{R}), \mathcal{D}_{\bar{\lambda},\mu})$ are computed by
Feigen and Fuchs \cite{ff}.
For $\bar{\lambda}\in \mathbb{R}^2$ the spaces $\mathrm{H}^1_\mathrm{diff}(\mathfrak{sl}(2),\mathcal{D}_{\bar{\lambda},\mu})$ are computed by
Bouarroudj \cite{b}. For $\bar{\lambda}\in \mathbb{R}^3$ the spaces $\mathrm{H}^1_\mathrm{diff}(\mathfrak{sl}(2),\mathcal{D}_{\bar{\lambda},\mu})$ are computed by O. Basdouri and N. Elamine \cite{be}. In this paper we are interested  to compute the  spaces
$\mathrm{H}^1_\mathrm{diff}(\mathfrak{sl}(2),\mathcal{D}_{\bar{\lambda},\mu})$ for $\bar{\lambda}\in \mathbb{R}^n$.

\section{Cohomology }

Let us first recall some fundamental concepts from cohomology theory~(see, e.g., \cite{Fu}). Let $\mathfrak{g}$ be a Lie algebra
acting on a space $V$ and let $\mathfrak{h}$ be a subalgebra of $\mathfrak{g}$. (If $\mathfrak{h}$ is omitted it assumed to be $\{0\}$). The space of
$\frak h$-relative $n$-cochains of $\mathfrak{g}$ with values in $V$ is the $\mathfrak{g}$-module
\begin{equation*}
C^n(\mathfrak{g},\mathfrak{h}, V ) := \mathrm{Hom}_{\frak
h}(\Lambda^n(\mathfrak{g}/\mathfrak{h}),V).
\end{equation*}
The {\it coboundary operator} $ \partial: C^n(\mathfrak{g},\mathfrak{h}, V)\longrightarrow C^{n+1}(\mathfrak{g},\mathfrak{h}, V )$ is a $\mathfrak{g}$-map
satisfying $\partial^2=0$. The operator $ \partial$  is defined by
$$\aligned
&(\partial f)(u_0,\dots,u_n)=\sum_{i=0}^n (-1)^iu_i f(u_0,\dots,\hat{\imath},\dots,u_n)+\\
&~~+ \sum_{0\leq i<j\leq n}(-1)^{i+j}f([u_i,u_j],u_0,\dots,\hat{\imath},\dots,\hat{\jmath},\dots, u_n).
\endaligned
$$
The kernel of $\partial|_{C^n}$, denoted $Z^n(\mathfrak{g},\mathfrak{h},V)$, is the space of $\frak h$-relative $n$-{\it cocycles}, among them,
the elements of $\partial(C^{n-1}(\mathfrak{g},\mathfrak{h}, V))$ are called $\mathfrak{h}$-relative $n$-{\it coboundaries}.
We denote $B^n(\mathfrak{g}, \mathfrak{h},V)$ the space of $n$-coboundaries.\\

By definition, the $n^{th}$ $\frak h$-relative  cohomology space is the quotient space
\begin{equation*}
\mathrm{H}^n (\mathfrak{g},\mathfrak{h},V)=Z^n(\mathfrak{g},\mathfrak{h},V)/B^n(\mathfrak{g},\mathfrak{h},V).
\end{equation*}
Here we consider $\mathrm{H}^n (\mathfrak{g},V)$ where $\mathfrak{g}=\mathfrak{sl}(2)$  and  $V=\mathcal{D}_{\bar{\lambda},\mu}$. In this paper we are interested
to the differential cohomology $\mathrm{H}^1_{\mathrm{diff}}$ (i.e., we consider only cochains that are given by differential operators).

\section{The spaces  $\mathrm{H}^1_\mathrm{diff}(\mathfrak{sl}(2),\mathcal{D}_{\bar{\lambda},\mu})$}
  Consider $\mu\in\mathbb{R}$, $\alpha=(\alpha_1\dots,\alpha_n)\in\mathbb{N}^n$ and $\bar{\lambda}=(\lambda_1\dots,\lambda_n)\in\mathbb{R}^n$, we define
  $$
  |\bar{\lambda}|=\sum_{i=1}^n\lambda_i,\quad \delta=\mu-|\bar{\lambda}|\quad\text{and}\quad|\alpha|=\sum_{i=1}^n\alpha_i.
  $$
  We also consider the canonical basis $(\varepsilon_1,\,\dots,\,\varepsilon_n)$ of $\mathbb{R}^n$ where $\varepsilon_i$ is the $n$-tuple
  $(0,\dots,0, 1,0,\dots, 0)$ ($1$ in the $i^{th}$ position).
 The space $\mathcal{D}_{\bar{\lambda},\mu}$ is spanned  by the operators $\Omega^\alpha$ defined by
$$
\Omega^\alpha(f_1dx^{\lambda_1}\otimes\cdots\otimes f_ndx^{\lambda_n})=f_1^{(\alpha_1)}\dots f_n^{(\alpha_n)}dx^{\mu}.
$$
\begin{lem}\label{lem0}
The action of $X_x$ on $\mathcal{D}_{\bar{\lambda},\mu}$ is diagonalizable. The operators $\Omega^\alpha$ are eigenvectors:
$$
X_x\cdot\Omega^\alpha=(\mu-|\bar{\lambda}|-|\alpha|)\Omega^\alpha=(\delta-|\alpha|)\Omega^\alpha.
$$
\end{lem}

The following lemma gives a first reduction of any 1-cocycle.
\begin{lem}\label{lem1} We have $$\mathrm{H}^1_\mathrm{diff}(\mathfrak{sl}(2),\mathcal{D}_{\bar{\lambda},\mu})=
\mathrm{H}^1_\mathrm{diff}(\mathfrak{sl}(2), X_1,\mathcal{D}_{\bar{\lambda},\mu}).$$
Moreover, up to a coboundary, any 1-cocycle $f\in\mathrm{Z}^1_\mathrm{diff}(\mathfrak{sl}(2),\mathcal{D}_{\bar{\lambda},\mu})$ can be expressed as follows:
\begin{equation*}
\label{cocycle}f(X_h) =\sum_\alpha B_\alpha h' \Omega^{\alpha} +\sum_\alpha C_\alpha h''\Omega^{\alpha},\qquad X_h\in\mathfrak{sl}(2),
\end{equation*}
where, for any $\alpha$, the coefficients $B_{\alpha+\varepsilon_i}$ and $C_\alpha$ are constants satisfying:
\begin{equation}\label{coef}
2(\delta-|\alpha|-1)C_\alpha+\sum_i(\alpha_i+1)(\alpha_i+2\lambda_i)B_{\alpha+\varepsilon_i}=0.
\end{equation}
\end{lem}
\begin{proofname}. Any 1-cocycle on $\mathfrak{sl}(2)$ should retains the following general form:
$$
f(X_h) = \sum_\alpha U_\alpha h \Omega^{\alpha}+\sum_\alpha V_\alpha h' \Omega^{\alpha}+\sum_\alpha W_\alpha h''\Omega^{\alpha},
$$
 where $U_\alpha$, $V_\alpha$ and $W_\alpha$ are, a priori, functions. First, we prove that the terms in $h$ can be annihilated by adding a coboundary.
Consider the $n$-ary differential operator
$$
g : \mathcal{F}_{\lambda_1}\otimes\cdots \otimes\mathcal{F}_{\lambda_n}\rightarrow\mathcal{F_\mu},\quad g=\sum_\alpha D_\alpha \Omega^{\alpha},\quad D_\alpha\in\mathbb{R}[x].
$$
 We have
\begin{equation}\label{delta}\begin{array}{lll}
\partial g(X_h) &=&hg'+\mu h'g-g\circ L^{\bar{\lambda}}_{X_h} \\[10pt]
&=& \displaystyle\sum_\alpha D'_\alpha h \Omega^{\alpha}+\displaystyle\sum_\alpha (\delta-|\alpha|)D_\alpha h' \Omega^{\alpha}\\[10pt]
&-&{1\over2}\displaystyle\sum_{\alpha}\sum_{i=1}^n\alpha_i(\alpha_i+2\lambda_i-1) D_{\alpha}h''\Omega^{\alpha-\varepsilon_i}. \end{array}\end{equation}
Thus, if $D'_\alpha=U_\alpha$ then $f-\partial g$ does not contain terms in h. So, we can replace $f$ by $f-\partial g$. That is, up to a coboundary,
any 1-cocycle on $\mathfrak{sl}(2)$ can be expressed as follows:
$$
f(X_h) = \sum_\alpha B_\alpha h' \Omega^{\alpha}+\sum_\alpha C_\alpha h''\Omega^{\alpha}.
$$

Now, for $X_{h_1},\,X_{h_2} \in\mathfrak{sl}(2)$, consider the 1-cocycle condition:
$$
f([X_{h_1}, X_{h_2}])- X_{h_1}\cdot f(X_{h_2})
+X_{h_2}\cdot f(X_{h_1}) = 0,
$$
which can be expressed as follows:
\begin{equation*}\begin{array}{lll} \displaystyle\sum_\alpha B'_\alpha ( h_1 h'_2-h'_1h_2) \Omega^{\alpha}+
\displaystyle\sum_\alpha C'_\alpha ( h_1 h''_2-h''_1h_2) \Omega^{\alpha}\\[10pt]
+{1\over2}\displaystyle\sum_\alpha\Big(2(\delta-|\alpha|-1)C_\alpha+\sum_{i=1}^n(\alpha_i+1)(\alpha_i+
2\lambda_i)B_{\alpha+\varepsilon_i}\Big)( h'_1 h''_2-h''_1h'_2)\Omega^{\alpha}=0.
\end{array}\end{equation*}
Thus, for all $\alpha$, we have $B'_\alpha=C'_\alpha=0$, and, moreover, the $B_{\alpha+\varepsilon_i}$ and $C_\alpha$ satisfy \eqref{coef}.

\hfill$\Box$
\end{proofname}\\

\begin{thm}\label{cor} \begin{itemize}
                         \item [1)] If $\delta\notin\mathbb{N}$ then ${\rm H}^1_\mathrm{diff}(\frak
{sl}(2),\mathcal{D}_{\bar{\lambda},\mu})=0$.
                         \item [2)] If $\delta=k\in\mathbb{N}$ then, up to a coboundary,  any 1-cocycle
$f\in\mathrm{Z}^1_\mathrm{diff}(\mathfrak{sl}(2),\mathcal{D}_{\bar{\lambda},\mu})$
can be expressed as follows:
\begin{equation}\label{cocyc}f(X_h) =\sum_{|\alpha|=k} B_\alpha h' \Omega^{\alpha}+\sum_{|\beta|=k-1} C_\beta h''\Omega^{\beta},
\qquad X_h\in\mathfrak{sl}(2),\end{equation}
where, for any $\beta$ such that $|\beta|=k-1$, the  $B_{\beta+\varepsilon_i}$  are constants satisfying:
\begin{equation}\label{coeff1}\sum_i(\beta_i+1)(\beta_i+2\lambda_i)B_{\beta+\varepsilon_i}=0.
\end{equation}
 \item [3)] If the $B_\alpha$ are not all zero then the cocycles \eqref{cocyc} are nontrivial.
                       \end{itemize}

\end{thm}

\begin{proofname}. Consider the 1-cocycle $f$ defined by (\ref{cocycle}) and consider the operator
$\partial g$ where
$$
g =\sum_{|\alpha|\neq \delta}{1\over \delta-|\alpha|} B_\alpha \Omega^{\alpha}.
$$
We have
\begin{equation*}\begin{array}{lll}\label{coboundary}
\partial g(X_h) &=&\displaystyle\sum_{|\alpha|\neq \delta} B_\alpha h' \Omega^{\alpha}-
{1\over2}\displaystyle\sum_{|\alpha|\neq \delta}\sum_{i=1}^n {\alpha_i(\alpha_i+2\lambda_i-1)\over \delta-|\alpha|}
B_\alpha h''\Omega^{\alpha-\varepsilon_i}\\ &=&\displaystyle\sum_{|\alpha|\neq \delta} B_\alpha h' \Omega^{\alpha}-
{1\over2}\displaystyle\sum_{|\beta|\neq \delta-1}\sum_{i=1}^n {(\beta_i+1)(\beta_i+2\lambda_i)\over \delta-|\beta|-1}
B_{\beta+\varepsilon_i} h''\Omega^{\beta}. \end{array}\end{equation*}
According to \eqref{coef} we have
\[
-{1\over2}\sum_{i=1}^n {(\beta_i+1)(\beta_i+2\lambda_i)\over \delta-|\beta|-1} B_{\beta+\varepsilon_i}=C_\beta.
\]
Therefore, we have
$$
\partial g(X_h)=\sum_{|\alpha|\neq \delta} B_\alpha h' \Omega^{\alpha}+\sum_{|\beta|\neq \delta-1} C_\beta h''\Omega^{\beta}
$$
and
$$
(f-\partial g)(X_h)=\sum_{|\alpha|=\delta} B_\alpha h' \Omega^{\alpha}+\sum_{|\beta|=\delta-1} C_\beta h''\Omega^{\beta}.
$$
Thus, if $\delta\notin\mathbb{N}$ then $f-\partial g=0$ (since $|\alpha|\in\mathbb{N}$ then $|\alpha|$ can not be equal to $\delta$), therefore ${\rm H}^1_\mathrm{diff}(\frak
{sl}(2),\mathcal{D}_{\bar{\lambda},\mu})=0$. If $\delta\in\mathbb{N}$ then the condition \eqref{coeff1} is coming from \eqref{coef}, since, in \eqref{cocyc} we have $\delta-|\beta|-1=0$.
Moreover, if $\delta=\mu-|\lambda|=|\alpha|=k$, then there are no terms in $h'\Omega^{\alpha}$ in the expression of $\partial g(X_h)$
for any $g\in \mathrm{D}_{\bar{\lambda}, \mu}$ (see \eqref{delta}), therefore, the non vanishing cocycle $f(X_h)=\sum_{|\alpha|=k} B_\alpha h' \Omega^{\alpha}$ are nontrivial.

\hfill$\Box$
\end{proofname}\\

Now, we prove that, generically, we can annihilate the term in $h''$ in the expression of the 1-cocycle \eqref{cocycle1} by adding a coboundary
and we describe completely the space $\mathrm{H}^1(\mathfrak{sl}(2),\mathcal{D}_{\bar{\lambda},\mu})$ for generic $\bar{\lambda}$.

\begin{thm}\label{third}
 If $\delta=k\in\mathbb{N}$ and $-2\bar{\lambda}\notin\{0,\,\dots,\,k-1\}^n$  then, we have
$$\mathrm{dim}\mathrm{H}^1_\mathrm{diff}(\mathfrak{sl}(2),\mathcal{D}_{\bar{\lambda},\mu})= \begin{pmatrix}n+k-2\\k\end{pmatrix}.$$
These spaces are spanned by the cocycles:
\begin{equation*}\label{cocycth}f(X_h) =\sum_{|\alpha|=k} B_\alpha h' \Omega^{\alpha},\end{equation*}
where the  $B_{\alpha}$  are constants satisfying the conditions \eqref{coeff1}.
   \end{thm}

\begin{proofname}.  By Theorem \ref{cor}, for $\delta=k\geq1$, any 1-cocycle $f$ can be expressed as follows:
\begin{equation}\label{cocycle1}f(X_h) =\sum_{|\alpha|=k} B_\alpha h' \Omega^{\alpha}+\sum_{|\beta|=k-1} C_\beta h''\Omega^{\beta}.\end{equation}
Consider
$$
g= \sum_{|\alpha|=k} D_\alpha \Omega^{\alpha}.
$$
We have
\begin{equation}\begin{array}{lll}\label{coboundary1}
\partial g(X_h) &=&-{1\over2}\displaystyle\sum_{|\alpha|=k}\sum_{i=1}^n\alpha_i(\alpha_i+2\lambda_i-1) D_{\alpha}h''\Omega^{\alpha-\varepsilon_i}\\ &=
&-{1\over2}\displaystyle\sum_{|\beta|=k-1}\sum_{i=1}^n(\beta_i+1)(\beta_i+2\lambda_i) D_{\beta+\varepsilon_i}h''\Omega^{\beta}. \end{array}
\end{equation}
Now, consider the linear system
 \begin{equation}\label{system}{1\over2}\displaystyle\sum_{i=1}^n(\beta_i+1)(\beta_i+2\lambda_i) D_{\beta+\varepsilon_i}=C_{\beta}, \quad |\beta|=k-1
 \end{equation}
 which express that
 $$
 \sum_{|\beta|=k-1} C_\beta h''\Omega^{\beta}=-\partial g(X_h).
 $$
 Without loss of generality, assume that $-2\lambda_1\notin\{0,\,\dots,\,k-1\}$. Choose arbitrarily the  $D_{\alpha}$ where
 $\alpha=(0,\alpha_2,\dots,\alpha_n)\in\mathbb{N}^n$ with $|\alpha|=k$. Now, for any $\alpha=(1,\alpha_2,\dots,\alpha_n)\in\mathbb{N}^n$ with
 $|\alpha|=k$, consider $\beta=(0,\alpha_2,\dots,\alpha_n)$ (then $|\beta|=k-1$). The coefficient $D_{\alpha}$ is uniquely defined by \eqref{system},
 in function of $C_\beta$ and the $D_{\beta+\varepsilon_i}$, for $i\geq2$, (indeed, $D_\alpha=D_{\beta+\varepsilon_1}$, $\beta_1=0$ and $\lambda_1\neq0$).
 Similarly, by \eqref{system}, we define the coefficients $D_{\alpha}$ with $\alpha_1=2$, in function of those with $\alpha_1=1$ and
 $C_{(1,\alpha_2,\dots,\alpha_n)}$. So, step by step, we define all the coefficients $D_{\alpha}$ so that $(f+\partial g)(X_h)$ does not contain terms
 in $h''$.\\

 The non vanishing 1-cocycles defined be \eqref{cocycth} are non trivial, since any coboundary does not contain terms in $h'$ (see \eqref{coboundary1}).
 Thus, the space $\mathrm{H}^1_\mathrm{diff}(\mathfrak{sl}(2),\mathcal{D}_{\bar{\lambda},\mu})$ is isomorphic to the space of solutions of the system of
 linear equations \eqref{coeff1}.  As before, we assume that $-2\lambda_1\notin\{0,\,\dots,\,k-1\}$, the space of solutions of the system of linear
 equations \eqref{coeff1} is generated by the arbitrary coefficients $B_{\alpha}$ where $\alpha=(0,\alpha_2,\dots,\alpha_n)\in\mathbb{N}^n$ with
 $|\alpha|=k$. The number of such $B_\alpha$ is the well known binomial coefficient with repetition $$\Gamma_{n-1}^k=\begin{pmatrix}n+k-2\\k\end{pmatrix},$$
 which is equal to the dimension of $\mathrm{H}^1_\mathrm{diff}(\mathfrak{sl}(2),\mathcal{D}_{\bar{\lambda},\mu}).$

\hfill$\Box$\end{proofname}

\section{Singular cases}

Now assume that
$$
-2\bar{\lambda}=(t_1,\,t_2,\,\dots,\,t_n)\in\{0,\,1,\,\dots,\,k-1\}^n.
$$
We cut the system of equations \eqref{coeff1} into the two following subsystems $(S_1)$ and $(S_{2})$:
\begin{equation}\label{s1}
(S_1):\quad\sum_{i=1}^n(\alpha_i+1)(\alpha_i+2\lambda_i)B_{\alpha+\varepsilon_i}=0,\quad \alpha_1\neq t_1.
\end{equation}
\begin{equation}\label{s2}
(S_2):\quad\sum_{i=1}^n(\alpha_i+1)(\alpha_i+2\lambda_i)B_{\alpha+\varepsilon_i}=0,\quad \alpha_1= t_1.
\end{equation}
From $(S_1)$ we extract the following system $(S'_1)$
\begin{equation}\label{s1prime}
(S'_1):\quad\sum_{i=2}^n(\alpha_i+1)(\alpha_i+2\lambda_i)B_{\alpha+\varepsilon_i}=0,\quad \alpha_1=t_1-1.
\end{equation}

The space $\mathrm{H}^1_\mathrm{diff}(\mathfrak{sl}(2),\mathcal{D}_{\bar{\lambda},\mu})$ is managed by the system \eqref{coeff1}, which is a system with $\Gamma_{n}^k$ unknowns and $\Gamma_{n}^{k-1}$ equations. Any equation is coming from a given $\alpha=(\alpha_1,\dots,\alpha_n)$ with $|\alpha|=k-1$. The rank of \eqref{coeff1} is then less or equal to $\Gamma_{n}^{k-1}$. The cocycles generating  $\mathrm{H}^1_\mathrm{diff}(\mathfrak{sl}(2),\mathcal{D}_{\bar{\lambda},\mu})$ are of two types:
$$
\text{type } 1:\quad\sum_{|\alpha|=k} B_\alpha h' \Omega^{\alpha}\quad\text{or}\quad\text{type } 2:\quad \sum_{|\beta|=k-1} C_\beta h'' \Omega^{\beta}.
$$
If the rank of \eqref{coeff1} is $\Gamma_{n}^{k-1}-\ell$ then the dimension of the space of cocycle of type 1 is
$$
\Gamma_{n}^k-\Gamma_{n}^{k-1}+\ell=\Gamma_{n-1}^{k}+\ell.
$$
The dimension of the space of classes of cocycles of type 2 is equal to the dimension of the space of all $C_\beta$, with $|\beta|=k-1$,
from which we subtract the rank of \eqref{system} (or \eqref{coeff1}).
Thus, the dimension of the space of classes of cocycles of type 2 is $\ell$ since the space of parameters $C_\beta$, with $|\beta|=k-1$,
is $\Gamma_{n}^{k-1}$-dimensional. Thus,
$$
\mathrm{dim}\mathrm{H}^1_\mathrm{diff}(\mathfrak{sl}(2),\mathcal{D}_{\bar{\lambda},\mu})=\Gamma_{n-1}^{k}+2\ell=\begin{pmatrix}n+k-2\\k\end{pmatrix}+2\ell.
$$

\begin{theorem}\label{th1} If $\sigma_n< k-1$, then
$$
\mathrm{dim}\mathrm{H}^1_\mathrm{diff}(\mathfrak{sl}(2),\mathcal{D}_{\bar{\lambda},\mu})=\begin{pmatrix}n+k-2\\k\end{pmatrix}.
$$

\end{theorem}
\begin{proofname}. This is equivalent to the fact that, in this case, the system \eqref{coeff1} is of maximal rank. We proceed by recurrence to prove this result. It is true for $n=2$ (see \cite{b}).
Assume that the result is true for $n-1$.  Consider the subsystems $(S_1)$ and $(S_2)$, it is easy to prove that $(S_1)$ is of maximal rank, while, by the recurrence hypothesis, the subsystem $(S_2)$ is also of maximal rank. Indeed, the subsystem $(S_2)$ can be considered as a system of equations related to the $(n-1)$-tuple $(\alpha_2,\dots,\alpha_n)$. So, we are in the case $n-1$ with
$$
\alpha_2+\cdots+\alpha_n=k-1-t_1=k'-1\quad\text{and}\quad \sigma_{n-1}=t_2+\cdots+t_n<k-t_1-1=k'-1.
$$
If $t_1=0$ then, combining the two subsystems $(S_1)$ and $(S_2)$, we get a system with maximal rank, since there are no common unknowns. Indeed, the unknowns $B_{(\alpha_1,\dots,\alpha_n)} $ of $(S_1)$ are all with $\alpha_1\neq0$, while those of $(S_2)$ are all with $\alpha_1=0$.

If $t_1>0$ then the unknowns $B_{(\alpha_1,\dots,\alpha_n)} $ of $(S_2)$ are all with $\alpha_1=t_1$. The unknowns $B_{(\alpha_1,\dots,\alpha_n)} $ with $\alpha_1=t_1$ can also appear in $(S_1)$, but they appear only in the equations corresponding to $(\alpha_1,\dots,\alpha_n) $ with $\alpha_1=t_1-1$. Therefore, we consider the system $(S'_1)$ which can be also considered as a system of equations related to $(\alpha_2,\dots,\alpha_n)$, but, here we have
 $$
 \alpha_2+\cdots+\alpha_n=k-t_1=k'-1\quad\text{and}\quad \sigma_{n-1}=t_2+\cdots+t_n<k-t_1-1=k'-2<k'-1.
 $$
For $(S'_1)$ we are in the case $n-1$, therefore, by the recurrence hypothesis, the system $(S'_1)$ is of maximal rank, therefore, there are no nontrivial combination of some equations of $(S_1)$ belonging to $(S_2)$, since there are no nontrivial combination of some equations of $(S_1)$ killing the $B_\alpha$ with $\alpha_1=t_1-1$. That is, the spaces of the equations respectively of $(S_1)$ and $(S_2)$ are supplementary. Thus, combining the two subsystems $(S_1)$ and $(S_2)$, we get a system of maximal rank which is the system \eqref{coeff1}. Theorem \ref{th1} is proved.

\hfill$\Box$
\end{proofname}
 \begin{theorem}\label{th2} If $\sigma_n=k-1$ then
$$
\mathrm{dim}\mathrm{H}^1_\mathrm{diff}(\mathfrak{sl}(2),\mathcal{D}_{\bar{\lambda},\mu})=\begin{pmatrix}n+k-2\\k\end{pmatrix}+2.
$$
\end{theorem}
\begin{proofname}. This is equivalent to the fact that, in this case, the system \eqref{coeff1} is of rank $\Gamma_{n}^{k-1}-1$. We proceed by recurrence to prove this fact. This is true for $n=2$ (see \cite{b}).
Assume that the result is true for $n-1$.  Consider the subsystems $(S_1)$ and $(S_2)$. As before $(S_1)$ is of maximal rank, it is of rank $\Gamma_{n}^{k-1}-\Gamma_{n-1}^{k-t_1-1}$.  The subsystem $(S_2)$ is of rank $\Gamma_{n-1}^{k-t_1-1}-1$ (by the recurrence hypothesis).

If $t_1=0$ then, as before, we see that the subsystems $(S_1)$ and $(S_2)$ are independent of each other. Therefore, combining the two subsystems $(S_1)$ and $(S_2)$, we get a system of rank $$\Gamma_{n}^{k-1}-\Gamma_{n-1}^{k-t_1-1}+\Gamma_{n-1}^{k-t_1-1}-1=\Gamma_{n}^{k-1}-1.$$ Indeed, the unknowns $B_{(\alpha_1,\dots,\alpha_n)} $ of $(S_1)$ are all with $\alpha_1\neq0$, while those of $(S_2)$ are all with $\alpha_1=0$ (there are common unknowns).

If $t_1>0$ then the unknowns $B_{(\alpha_1,\dots,\alpha_n)} $ of $(S_2)$ are all with $\alpha_1=t_1$. The unknowns $B_{(\alpha_1,\dots,\alpha_n)} $ with $\alpha_1=t_1$ appear in $(S_1)$ only in the equations corresponding to ${(\alpha_1,\dots,\alpha_n)} $ with $\alpha_1=t_1-1$. As before, we consider the system $(S'_1)$. We are in the case $n-1$ with $$\alpha_2+\cdots+\alpha_n=k-t_1=k'-1\quad\text{and}\quad \sigma_{n-1}=t_2+\cdots+t_n=k-t_1-1=k'-2<k'-1.$$

By the recurrence hypothesis, the system $(S'_1)$ is of maximal rank, therefore, there are no nontrivial combination of some equations of $(S_1)$ killing the $B_\alpha$ with $\alpha_1=t_1-1$, therefore, combining the two subsystems $(S_1)$ and $(S_2)$, we get a system with rank
$$\Gamma_{n}^{k-1}-\Gamma_{n-1}^{k-t_1-1}+\Gamma_{n-1}^{k-t_1-1}-1=\Gamma_{n}^{k-1}-1.$$
Theorem \ref{th2} is proved.

\hfill$\Box$
\end{proofname}
\begin{remark} If $\sigma_n=k-1$, we see that the equation corresponding to ${(t_1,\dots,t_n)} $ is trivial, so, it's enough to prove that, if we subtract this equation from \eqref{coeff1}, we get a maximal rank system. In this case the space of cocycles of  type 2 is one dimensional, spanned by
\begin{equation*} \omega(X_h) =h''\Omega^{(t_1,\dots, t_n)}.\end{equation*}
\end{remark}

\begin{theorem}\label{th3} If $\sigma_n= k$ then
$$
\mathrm{dim}\mathrm{H}^1_\mathrm{diff}(\mathfrak{sl}(2),\mathcal{D}_{\bar{\lambda},\mu})=\begin{pmatrix}n+k-2\\k\end{pmatrix}+2(s-1),
$$
where $s$ is the number of  $t_i\geq1$.
\end{theorem}

\begin{proofname}. We proceed by recurrence to prove that the rank of \eqref{coeff1} is
$$
\Gamma_{n}^{k-1}-(s-1)=\begin{pmatrix}n+k-2\\k-1\end{pmatrix}-(s-1).
$$
This is true for $n=2$ (see \cite{b}), indeed, for $n=2$ we have necessarily $s=2$.
Assume that the result is true for $n-1$. As before, if $t_1=0$ then the subsystems $(S_1)$ and $(S_2)$ are independent of each other. The subsystem $(S_1)$ is of rank $\Gamma_{n}^{k-1}-\Gamma_{n-1}^{k-1}$, while, according to the recurrence hypothesis, the subsystem $(S_2)$ is of rank $\Gamma_{n-1}^{k-1}-(s-1).$
Therefore, the rank of \eqref{coeff1} is
$$
\Gamma_{n}^{k-1}-\Gamma_{n-1}^{k-1}+\Gamma_{n-1}^{k-1}-(s-1)=\Gamma_{n}^{k-1}-(s-1)=\begin{pmatrix}n+k-2\\k-1\end{pmatrix}-(s-1).
$$

Now, for $t_1>0$, the subsystem $(S_1)$ is of rank $\Gamma_{n}^{k-1}-\Gamma_{n-1}^{k-t_1-1}$, while, according to the recurrence hypothesis, the subsystem $(S_2)$ is of rank
$$
\Gamma_{n-1}^{k-t_1-1}-(s-2).
$$
Consider the subsystem $(S'_1)$. We are in the case $n-1$ with
$$
\alpha_2+\cdots+\alpha_n=k-t_1=k'-1\quad\text{and}\quad \sigma_{n-1}=t_2+\cdots+t_n=k-t_1=k'-1.
$$
Therefore, the subsystem $(S'_1)$ is of rank
$$
\Gamma_{n-1}^{k-t_1}-1.
$$
The equation corresponding to ${(t_1-1,t_2,\dots,t_n)}$ appear as a trivial equation in $(S'_1)$, corresponding in $(S_1)$ to the equation
\begin{equation}\label{last}
B_{(t_1,t_2,\dots,t_n)}=0.
\end{equation}
The equation \eqref{last} appears also in $(S_2)$ corresponding to any ${(t_1,t_2,\dots,t_n)-\varepsilon_i}$ for $i\geq2$.

Obviously, if we subtract the trivial equation from $(S'_1)$ we obtain a system of maximal rank.
Thus, the system \eqref{coeff1} is of rank
$$
\Gamma_{n}^{k-1}-\Gamma_{n-1}^{k-t_1-1}+\Gamma_{n-1}^{k-t_1-1}-(s-2)-1=\Gamma_{n}^{k-1}-(s-1).
$$
Theorem \ref{th3} is proved.

\hfill$\Box$
\end{proofname}
\begin{remark} The result of the previous theorem can be explained by the fact that the equations \eqref{coeff1} corresponding to $(t_1,\dots,t_n)-\varepsilon_i$, for $t_i>0$, are the same. They are all equivalent to \eqref{last}. So, it's enough to prove that, if we subtract $(s-1)$ equations from these $s$ equivalent equations we get a maximal rank system. The space of cocycles of type 2 is spanned by
\begin{equation*} \Gamma_i(X_h) =h''\Omega^{\gamma_i},\quad\gamma_i=(t_1,t_2,\dots,t_n)-\varepsilon_i,\quad t_i\geq1.\end{equation*}
\end{remark}

\begin{theorem}\label{th4} If $\sigma_n= k+1$, then
$$
\mathrm{dim}\mathrm{H}^1_\mathrm{diff}(\mathfrak{sl}(2),\mathcal{D}_{\bar{\lambda},\mu})=\left\{\begin{array}{ll}\Gamma_{n}^{k-1}+s(s-1)-2r
\quad&\text{if}\quad\max{t_i}\geq2\\[8pt]
\Gamma_{n}^{k-1}\quad&\text{if}\quad\max{t_i}=1.
\end{array}\right.
$$
where $s$ is the number of $t_i\geq1$ and $r$ is the number of $t_i=1$.
\end{theorem}

\begin{proofname}. Assume that $t_1\geq1$ and consider the system $(S'_1)$. We are in the case $n-1$ with
$$
\alpha_2+\cdots+\alpha_n=k-t_1=k'-1\quad\text{and}\quad \sigma_{n-1}=t_2+\cdots+t_n=k+1-t_1=k'.
$$
Therefore, it was proved in the proof of Theorem \ref{th3} that the system $(S'_1)$ is of rank
$$
\Gamma_{n-1}^{k-t_1}-(s-2).
$$
In $(S'_1)$ the equations corresponding to ${(t_1-1,t_2,\dots,t_n)-\varepsilon_{i}}$, for $t_i\geq1$, are equivalent to
$$
B_{(t_1-1,t_2,\dots,t_n)}=0.
$$
But, the correspondent equations in $(S_1)$ are
\begin{equation}\label{lastt}
B_{(t_1,t_2,\dots,t_n)-\varepsilon_{i}}+B_{(t_1-1,t_2,\dots,t_n)}=0,\quad t_i\geq1.
\end{equation}

{\bf Case 1: $\max t_i\geq2$.}  Assume that $t_2\geq2$. We proceed by recurrence to prove that the rank of \eqref{coeff1} is
$$
\Gamma_{n}^{k-1}-\frac{1}{2}s(s-1)+r.
$$
This is true for $n=2$ (see \cite{b}), indeed, for $n=2$ we have necessarily $s=2$ and $r=0$. Assume that the result is true for $n-1$.

\bigskip

 Assume that $t_1=1$. By the recurrence hypothesis, the subsystem $(S_2)$ is of rank
$$
\Gamma_{n-1}^{k-t_1-1}-\frac{1}{2}(s-2)(s-1)+r-1.
$$
We have $t_2\geq2$, then, in $(S_2)$ the equation corresponding to ${(t_1,t_2-2,\dots,t_n)}$ gives
$$
B_{(t_1,t_2-1,\dots,t_n)}=0.
$$
Therefore, the equation corresponding to ${(t_1,t_2-1,\dots,t_n)-\varepsilon_{i}}$, for $t_i\geq1$, gives
\begin{equation}\label{lasst}
B_{(t_1,t_2,\dots,t_n)-\varepsilon_{i}}=0,\quad t_i\geq1,
\end{equation}
and the $(s-2)$ correspondent equations \eqref{lastt} in $(S_1)$, for $i\geq3$, become trivial (we can also say that the equations \eqref{lastt} of $(S_1)$, for $i\geq3$, are combination of some equations of $(S_2)$).
If we subtract the equations \eqref{lasst}, for $i\geq3$, from $(S'_1)$,  we obtain a maximal rank system.
Therefore, the rank of \eqref{coeff1} is
$$
\Gamma_{n}^{k-1}-\Gamma_{n-1}^{k-t_1-1}+\Gamma_{n-1}^{k-t_1-1}-\frac{1}{2}(s-2)(s-1)+r-1-(s-2)=\begin{pmatrix}n+k-2\\k\end{pmatrix}-\frac{1}{2}s(s-1)+r.
$$

If $t_1\geq2$, then the subsystem $(S_2)$ is of rank
$$
\Gamma_{n-1}^{k-t_1-1}-\frac{1}{2}(s-2)(s-1)+r.
$$
But, here we have $B_{(t_1-1,t_2,\dots,t_n)}=0$ as equation corresponding to $(t_1-2,t_2,\dots,t_n)$. Therefore, the equation corresponding to $(t_1-1,t_2,\dots,t_n)-\varepsilon_i$ gives $B_{(t_1,t_2,\dots,t_n)-\varepsilon_i}=0$. So, the $(s-2)$ correspondent equations \eqref{lastt} in $(S_1)$ (for $i\geq3$) become trivial, but, the equation in $(S_1)$: $B_{(t_1,t_2-1,\dots,t_n)}=0$, appear also in $(S_2)$ as equation corresponding to ${(t_1,t_2-2,\dots,t_n)}$. Thus, the rank of \eqref{coeff1} is
$$
\Gamma_{n}^{k-1}-\Gamma_{n-1}^{k-t_1-1}+\Gamma_{n-1}^{k-t_1-1}-\frac{1}{2}(s-2)(s-1)+r-(s-2)-1=\begin{pmatrix}n+k-2\\k\end{pmatrix}-\frac{1}{2}s(s-1)+r.
$$

{\bf Case 2: $\max t_i=1$.} In this case we prove the rank of \eqref{coeff1} is
$$
\Gamma_{n}^{k-1}=\begin{pmatrix}n+k-2\\k\end{pmatrix}.
$$
 By the recurrence hypothesis, the subsystem $(S_2)$ is of rank
$$
\Gamma_{n-1}^{k-t_1-1}.
$$

 Assume that $t_1=t_2=1$. Therefore, the $(s-2)$ equations, in $(S'_1)$, corresponding to ${(t_1-1,t_2,\dots,t_n)-\varepsilon_{i}}$, $t_i=1$ for $i\geq3$, are equivalent to $B_{(t_1-1,t_2,\dots,t_n)}=0$. Obviously, if we subtract these equations from $(S'_1)$, we get a maximal rank system.
Therefore, the rank of \eqref{coeff1} is
$$
\Gamma_{n}^{k-1}-\Gamma_{n-1}^{k-t_1-1}+\Gamma_{n-1}^{k-t_1-1}=\begin{pmatrix}n+k-2\\k\end{pmatrix}.
$$
Theorem \ref{th4} is proved.

\hfill$\Box$
\end{proofname}

\bigskip

\begin{remark} In the previous theorem, if there exist some $t_i\geq2$, then, for any $t_{i}\geq2$, \eqref{coeff1} applied to ${(t_1,t_2,\dots,t_n)-2\varepsilon_{i}}$ gives  $B_{(t_1,t_2,\dots,t_n)-\varepsilon_{i}}=0$. We use one of $t_i\geq2$ to prove that $B_{(t_1,t_2,\dots,t_n)-\varepsilon_{j}}=0$ for $t_j=1$. That is, we
have $r$ equations: $B_{(t_1,t_2,\dots,t_n)-\varepsilon_{j}}=0$,
corresponding to ${(t_1,t_2,\dots,t_n)-\varepsilon_{i_0}-\varepsilon_j}$ for a fixed $t_{i_0}\geq2$ and $t_j=1$.
Therefore, all other equations corresponding to ${(t_1,t_2,\dots,t_n)-\varepsilon_{i}-\varepsilon_j}$, for $i< j$ and $t_i,\,t_j\geq1$, become trivial. We prove that if we subtract these trivial equations we get a maximal rank system. Thus, the rank of \eqref{coeff1} is
$$
\Gamma_{n}^{k-1}-\frac{1}{2}s(s-1)+r=\begin{pmatrix}n+k-2\\k\end{pmatrix}-\frac{1}{2}s(s-1)+r.
$$
\end{remark}

\begin{theorem}\label{th5} If $\sigma_n= k+2$, then
$$
\mathrm{dim}\mathrm{H}^1_\mathrm{diff}(\mathfrak{sl}(2),\mathcal{D}_{\bar{\lambda},\mu})=\begin{pmatrix}n+k-2\\k\end{pmatrix}+s(s-1),
$$
where $s$ is the number of $t_i\geq3$.
\end{theorem}

\begin{proofname}. We proceed by recurrence to prove that the rank of \eqref{coeff1} is
$$
\Gamma_{n}^{k-1}-\frac{1}{2}s(s-1)=\begin{pmatrix}n+k-2\\k\end{pmatrix}-\frac{1}{2}s(s-1).
$$
This is true for $n=2$. Assume that it is true for $n-1$. As usually, the result is true if $t_1=0$.

\bigskip

If $\max t_i=1$ then  $(S_2)$ is of rank
$$
\Gamma_{n-1}^{k-t_1-1},
$$
and the system $(S'_1)$ is of maximal rank. Therefore, the rank of \eqref{coeff1} is
$$
\Gamma_{n}^{k-1}-\Gamma_{n-1}^{k-t_1-1}+\Gamma_{n-1}^{k-t_1-1}=\begin{pmatrix}n+k-2\\k\end{pmatrix}.
$$

\bigskip

If $\max t_i=2$ then  we can assume that $t_1=2$. In this case $(S_2)$ is of rank
$$
\Gamma_{n-1}^{k-t_1-1}.
$$
For $(S'_1)$ we distinguish two cases. $\max_{i>1}t_i=1$ or $\max_{i>1}t_i=2$. In the first case, the system $(S'_1)$ is of maximal rank. Therefore, the rank of \eqref{coeff1} is
$$
\Gamma_{n}^{k-1}-\Gamma_{n-1}^{k-t_1-1}+\Gamma_{n-1}^{k-t_1-1}=\begin{pmatrix}n+k-2\\k\end{pmatrix}.
$$
In the second case, the system $(S'_1)$ is of rank
$$
\Gamma_{n-1}^{k-t_1}-\frac{1}{2}(s'-2)(s'-1)+r',
$$
where $s'$ is the number of $t_i\geq1$ and $r'$ is the number of $t_i=1$. Indeed, for the system $(S'_1)$, we are in the case $n-1$ with
$$
\alpha_2+\cdots+\alpha_n=k-t_1=k'-1\quad\text{and}\quad \sigma_{n-1}=t_2+\cdots+t_n=k+2-t_1=k'+1.
$$
Assume that $t_2=2$. In the system $(S'_1)$, for $i\geq2$ such that $t_i\geq2$, the equations corresponding to  ${(t_1-1,t_2,\dots,t_n)-2\varepsilon_{i}}$ are
\begin{equation}\label{last3}
B_{(t_1-1,t_2,\dots,t_n)-\varepsilon_{i}}=0,\quad t_i\geq2,\quad i\geq2.
\end{equation}
In particular, we have $B_{(t_1-1,t_2-1,\dots,t_n)}=0.$
Therefore, there are $r'$ equations
\begin{equation}\label{last4}
B_{(t_1-1,t_2,\dots,t_n)-\varepsilon_{i}}=0,\quad t_i=1,\quad i\geq3,
\end{equation}
corresponding to ${(t_1-1,t_2-1,\dots,t_n)-\varepsilon_{i}}$, for $i\geq3$ and $t_i=1$. Therefore, all other equations corresponding to ${(t_1-1,t_2,\dots,t_n)-\varepsilon_{i}-\varepsilon_j}$, for $i< j$ and $t_i,\,t_j\geq1$, become trivial in $(S'_1)$. The number of these trivial equations is  $\frac{1}{2}(s'-2)(s'-1)-r'$. Of course, if we subtract these trivial equations from $(S'_1)$ we obtain a maximal rank system.
Therefore, the rank of \eqref{coeff1} is
$$
\Gamma_{n}^{k-1}-\Gamma_{n-1}^{k-t_1-1}+\Gamma_{n-1}^{k-t_1-1}=\begin{pmatrix}n+k-2\\k\end{pmatrix}.
$$
Indeed, these trivial equations of $(S'_1)$ are associated in $(S_1)$ to
\begin{equation}\label{tt1}
B_{(t_1,t_2,\dots,t_n)-\varepsilon_{i}-\varepsilon_j}=0.
\end{equation}
But, \eqref{tt1} appears also in $(S_2)$ as equation corresponding to ${(t_1,t_2-1,\dots,t_n)-\varepsilon_{i}-\varepsilon_j}$, since we have also $B_{(t_1,t_2-1,\dots,t_n)-\varepsilon_{i}}=0$, for $t_i\geq2$.

\bigskip

Now, if $t_1\geq3$, then $(S_2)$ is of rank
$$
\Gamma_{n-1}^{k-t_1-1}-\frac{1}{2}(s-2)(s-1).
$$
The rank of $(S'_1)$ is
$$
\Gamma_{n-1}^{k-t_1}-\frac{1}{2}(s'-2)(s'-1)+r'.
$$
In $(S'_1)$, for $i\geq2$ such that $t_i\geq2$, the equations corresponding to  ${(t_1-1,t_2,\dots,t_n)-2\varepsilon_{i}}$ are
\begin{equation}\label{last5}
B_{(t_1-1,t_2,\dots,t_n)-\varepsilon_{i}}=0,\quad t_i\geq2,\quad i\geq2.
\end{equation}
Moreover, there are $r'$ equations
\begin{equation}\label{last6}
B_{(t_1-1,t_2,\dots,t_n)-\varepsilon_{i}}=0,\quad t_i=1,\quad i\geq2,
\end{equation}
corresponding to ${(t_1-1,t_2,\dots,t_n)-\varepsilon_{i_0}-\varepsilon_{i}}$, for a fixed $i_0\geq2$ such that $t_{i_0}\geq2$ and $i\geq2$ such that $t_i=1$. All other equations corresponding to ${(t_1-1,t_2,\dots,t_n)-\varepsilon_{i}-\varepsilon_j}$, for $i< j$ such that $t_i,\,t_j\geq1$, become trivial in $(S'_1)$. The number of these trivial equations is  $\frac{1}{2}(s'-2)(s'-1)-r'$. These trivial equations in $(S'_1)$ appear in $(S_1)$ as
\begin{equation}\label{ss2}
B_{(t_1,t_2,\dots,t_n)-\varepsilon_{i}-\varepsilon_j}=0.
\end{equation}

For any $t_{i}\geq3$, the equation corresponding to $(t_1,t_2,\dots,t_n)-3\varepsilon_{i}$ gives  $B_{(t_1,t_2,\dots,t_n)-2\varepsilon_{i}}=0$. Therefore, the equation \eqref{ss2} appear in $(S_2)$ as equation corresponding to $(t_1,t_2,\dots,t_n)-3\varepsilon_{i}$ only for $i\geq2$ such that $t_i\geq3$. Thus, the rank of \eqref{coeff1} is
$$
\Gamma_{n}^{k-1}-\Gamma_{n-1}^{k-t_1-1}+\Gamma_{n-1}^{k-t_1-1}-\frac{1}{2}(s-2)(s-1)-(s-1)=\begin{pmatrix}n+k-2\\k\end{pmatrix}-\frac{1}{2}s(s-1).
$$
Theorem \ref{th5} is proved.

\hfill$\Box$
\end{proofname}

\begin{theorem}\label{th6} If $\sigma_n= k+m$, with $m\geq2$, then
$$
\mathrm{dim}\mathrm{H}^1_\mathrm{diff}(\mathfrak{sl}(2),\mathcal{D}_{\bar{\lambda},\mu})=\begin{pmatrix}n+k-2\\k\end{pmatrix}+s(s-1),
$$
where $s$ is the number of $t_i> m$.
\end{theorem}

\begin{proofname}.  Assume that $m\geq3$, since the case $m=2$ was treated in the previous theorem. We proceed by recurrence to prove that the rank of \eqref{coeff1} is
$$
\Gamma_{n}^{k-1}-\frac{1}{2}s(s-1).
$$
This is true for $n=2$. Assume that it is true for $n-1$.
\subsection{If $t_1\leq m$.} In this case the system $(S_2)$ is of rank
$$
\Gamma_{n-1}^{k-t_1-1}-\frac{1}{2}s(s-1).
$$
The system $(S'_1)$ is of rank
$$
\Gamma_{n-1}^{k-t_1}-\frac{1}{2}s'(s'-1)
$$
where $s'$ is the number of $t_i\geq m$ for $i\geq2$. Indeed, for the system $(S'_1)$, we are in the case $n-1$ with
$$
\alpha_2+\cdots+\alpha_n=k-t_1=k'-1\quad\text{and}\quad \sigma_{n-1}=t_2+\cdots+t_n=k+m-t_1=k'+m-1.
$$
In $(S'_1)$, for any $t_{i}\geq m$, the equation corresponding to $\alpha=(t_1-1,t_2,\dots,t_n)-m\varepsilon_{i}$ is
\begin{equation}\label{mmm}
B_{(t_1-1,t_2,\dots,t_n)-(m-1)\varepsilon_{i}}=0.
\end{equation}

\bigskip

{\bf Case 1:} $m=2h+1$ with $h\geq1$.
 In $(S'_1)$, for $i\neq j$, $t_i\geq 2h+1$ and $t_j\geq1$, according to \eqref{mmm}, the equation corresponding to $(t_1-1,t_2,\dots,t_n)-2h\varepsilon_{i}-\varepsilon_j$ gives
 $$
 B_{(t_1-1,t_2,\dots,t_n)-(2h-1)\varepsilon_{i}-\varepsilon_j}=0.
 $$
 Step by step, for $t_i\geq(2h+1)$ and $t_j\geq(2h+1)$, $i\neq j$, the equations corresponding respectively to   $(t_1-1,t_2,\dots,t_n)-(h+1)\varepsilon_{i}-h\varepsilon_j$ and $(t_1-1,t_2,\dots,t_n)-h\varepsilon_{i}-(h+1)\varepsilon_j$ are equivalent to the same equation which is
\begin{equation}\label{ab1}
B_{(t_1-1,t_2,\dots,t_n)-h\varepsilon_{i}-h\varepsilon_j}=0.
\end{equation}
Therefore, if we subtract from $(S'_1)$ the $\frac{1}{2}s'(s'-1)$ equations corresponding to  $(t_1-1,t_2,\dots,t_n)-(h+1)\varepsilon_{i}-h\varepsilon_j$, where $t_i,\,t_j\geq(2h+1)$ and $i< j$, we get a maximal rank. Thus, the rank of \eqref{coeff1} is
$$
\Gamma_{n}^{k-1}-\Gamma_{n-1}^{k-t_1-1}+\Gamma_{n-1}^{k-t_1-1}-\frac{1}{2}s(s-1)=\begin{pmatrix}n+k-2\\k\end{pmatrix}-\frac{1}{2}s(s-1).
$$
\bigskip

{\bf Case 2: $ m=2h$}.
 In $(S'_1)$, for $i\neq j$, $t_i\geq m$ and $t_j\geq1$, according to \eqref{mmm}, the equation corresponding to $(t_1-1,t_2,\dots,t_n)-(2h-1)\varepsilon_{i}-\varepsilon_j$ gives
 $$
 B_{(t_1-1,t_2,\dots,t_n)-(2h-2)\varepsilon_{i}-\varepsilon_j}=0.
 $$
 Consider a fixed $t_{i_0}>1$. The equation corresponding to $(t_1-1,t_2,\dots,t_n)-(2h-2)\varepsilon_{i}-\varepsilon_{i_0}-\varepsilon_j$, for $t_i\geq m$, gives
 $$
 B_{(t_1-1,t_2,\dots,t_n)-(2h-3)\varepsilon_{i}-\varepsilon_{i_0}-\varepsilon_j}=0.
 $$
 Step by step, for $t_i,\,t_j\geq m$, $i\neq j$, the equations corresponding respectively to   $(t_1-1,t_2,\dots,t_n)-\varepsilon_{i_0}-(h-1)\varepsilon_{i}-h\varepsilon_j$ and $(t_1-1,t_2,\dots,t_n)-\varepsilon_{i_0}-h\varepsilon_{i}-(h-1)\varepsilon_j$ are equivalent to the same equation which is
\begin{equation}\label{ab2}
B_{(t_1-1,t_2,\dots,t_n)-\varepsilon_{i_0}-(h-1)\varepsilon_{i}-(h-1)\varepsilon_j}=0.
\end{equation}
Thus, we have the same result as in the previous case.
\subsection{If $t_1>m$} In this case the system $(S_2)$ is of rank
$$
\Gamma_{n-1}^{k-t_1-1}-\frac{1}{2}(s-2)(s-1).
$$
The system $(S'_1)$ is of rank
$$
\Gamma_{n-1}^{k-t_1}-\frac{1}{2}s'(s'-1)
$$
where $s'$ is the number of $t_i\geq m$ for $i\geq2$. We proceed as in the previous case, but here, since $t_1>m$, we prove that in $(S'_1)$ the equations \eqref{ab1} and \eqref{ab2} become trivial for $i=j$ and $t_i>m$. The corresponding $(s-1)$ equations in $(S_1)$ are respectively
$$
B_{(t_1,t_2,\dots,t_n)-(2h+1)\varepsilon_{i}}=0\quad\text{and}\quad B_{(t_1,t_2,\dots,t_n)-\varepsilon_{i_0}-(2h-1)\varepsilon_{i}}=0.
$$
But, these equations appear also in $(S_2)$ as equations corresponding respectively to $(t_1,t_2,\dots,t_n)-(2h+2)\varepsilon_{i}$ and $(t_1,t_2,\dots,t_n)-\varepsilon_{i_0}-2h\varepsilon_{i}$ for any $i\geq2$ such that $t_i>m$. Thus, the rank of \eqref{coeff1} is
$$
\Gamma_{n}^{k-1}-\Gamma_{n-1}^{k-t_1-1}+\Gamma_{n-1}^{k-t_1-1}-\frac{1}{2}(s-2)(s-1)-(s-1)=\begin{pmatrix}n+k-2\\k\end{pmatrix}-\frac{1}{2}s(s-1).
$$
Theorem \ref{th6} is proved

\hfill$\Box$
\end{proofname}

\bigskip

Note that for $n=2$ and $\sigma_2\geq k-1$, we have always $$s(s-1)=2(s-1)=s(s-1)-2r=2$$ and we know that, in this case, we have
$\mathrm{H}^1_\mathrm{diff}(\mathfrak{sl}(2),\mathcal{D}_{\bar{\lambda},\mu})=3=\begin{pmatrix}2+k-2\\k\end{pmatrix}+2$.

\bigskip

Now, we summarize our results in the following theorem
\begin{theorem}   $$
\mathrm{dim}\mathrm{H}^1_\mathrm{diff}(\mathfrak{sl}(2),\mathcal{D}_{\bar{\lambda},\mu})=\left\{\begin{array}{lllll}\Gamma_{n-1}^k\quad&\text{if}\quad\sigma_n <k-1,\\[6pt]
\Gamma_{n-1}^k+2\quad&\text{if}\quad\sigma_n=k-1,\\[6pt]
\Gamma_{n-1}^k+2(s-1)\quad&\text{if}\quad\sigma_n=k,\\[6pt]
\Gamma_{n-1}^k+(s+r)(s+r-1)-2r\quad&\text{if}\quad\sigma_n=k+1,\,\text{ and }\,\max{t_i}\geq2\\[6pt]
\Gamma_{n-1}^k\quad&\text{if}\quad\sigma_n=k+1,\,\text{ and }\,\max{t_i}=1\\[6pt]
\Gamma_{n-1}^k+s(s-1)\quad&\text{if}\quad\sigma_n=k+m,\quad m\geq2,\\[6pt]
\end{array}\right..
$$
where $s$ is the number of $t_i> \sigma_n-k$, $r$ is the number of $t_i=1$ and $\Gamma_{n-1}^k=\begin{pmatrix}n+k-2\\k\end{pmatrix}$.
\end{theorem}


\end{document}